\documentclass[12pt,reqno]{amsart}
\usepackage{amsfonts}
\usepackage{amsthm}
\usepackage{amsmath,amssymb}
\usepackage{amscd}
\usepackage[latin2]{inputenc}
\usepackage{t1enc}
\usepackage[mathscr]{eucal}
\usepackage{indentfirst}
\usepackage{graphicx}
\usepackage{graphics}
\usepackage{pict2e}
\usepackage{epic}
\usepackage[margin=2.9cm]{geometry}
\usepackage{epstopdf} 
\usepackage{bm}
\usepackage[all]{xy}
\usepackage{mathrsfs}
\usepackage{lmodern}
\newtheoremstyle{break}%
{}{}%
{\itshape}{}%
{\bfseries}{}
{\newline}{}

\newtheoremstyle{dbreak}%
{}{}%
{\upshape}{}%
{\bfseries}{}
{\newline}{}

\theoremstyle{break}
\newtheorem{Thm}{Theorem}[section]
\newtheorem{Lem}[Thm]{Lemma}
\newtheorem{Prop}[Thm]{Proposition}
\newtheorem{Cor}[Thm]{Corollary}

\theoremstyle{dbreak}
\newtheorem{Def}[Thm]{Definition}
\newtheorem{Rem}[Thm]{Remark}

\newcommand{\prf}{\noindent\underline{$Proof.$}}

\numberwithin{equation}{section}

\newenvironment{nouppercase}{%
\renewcommand{\uppercasenonmath}[1]{}}{}

\makeatletter
\@namedef{subjclassname@2020}{%
  \textup{2020} Mathematics Subject Classification}
\makeatother

\begin{document}

\title[Geometric properties of homomorphisms between the absolute Galois groups]{Geometric properties of homomorphisms between the absolute Galois groups of mixed-characteristic complete discrete valuation fields with perfect residue fields}
\author{Takahiro Murotani}



\address[Takahiro Murotani]{Department of Mathematics, Tokyo Institute of Technology, Tokyo 152-8551, Japan}

\email{murotani.t.aa@m.titech.ac.jp}

\begin{abstract}
Although the analogue of the theorem of Neukirch-Uchida for $p$-adic local fields fails to hold as it is,  Mochizuki proved a certain analogue of this theorem for the absolute Galois groups with ramification filtrations of $p$-adic local fields.
Moreover, Mochizuki and Hoshi gave various (necessary and) sufficient conditions for homomorphisms between the absolute Galois groups of $p$-adic local fields to be ``geometric'' (i.e., to arise from homomorphisms of fields).
In the present paper, we consider similar problems for general mixed-characteristic complete discrete valuation fields with perfect residue fields.
One main result gives (necessary and) sufficient conditions for homomorphisms between the absolute Galois groups of mixed-characteristic complete discrete valuation fields with residue fields algebraic over the prime fields to be geometric.
We also give a ``weak-Isom'' anabelian result for homomorphisms between the absolute Galois groups of these fields.
\end{abstract}

\begin{nouppercase}
\maketitle
\end{nouppercase}

\tableofcontents

\section*{Introduction}

In anabelian geometry, we consider how much information of the geometry of schemes is recovered from their \'{e}tale fundamental groups.
Even in the case where the schemes are the spectra of fields (so, in this case, we consider the absolute Galois groups of the fields), this problem is highly non-trivial.
For number fields, an affirmative result, the theorem of Neukirch-Uchida (cf., e.g., \cite[Corollary 12.2.2]{NSW}), is known.
(Moreover, by using this result, Hoshi reconstructed number fields from their absolute Galois groups in the sense of mono-anabelian reconstruction (cf. \cite[Theorem A]{mono}).)
More generally, a similar result for infinite fields finitely generated over prime fields was proved by Pop (cf. \cite[Theorem 2]{P}).
Grothendieck, who is the originator of anabelian geometry, considered that anabelian geometry should be developed over fields finitely generated over prime fields.
So, according to his original conjecture, it seems that these results completed anabelian geometry of fields.
Moreover, it is difficult to generalize these results to other kinds of fields since, for example, the analogue of the theorem of Neukirch-Uchida for $p$-adic local fields fails to hold as it is.
However, Mochizuki proved a certain analogue of the theorem of Neukirch-Uchida for the absolute Galois groups with ramification filtrations of $p$-adic local fields (cf. \cite[Theorem 4,2]{version}).
(Moreover, Abrashkin proved an analogue of the theorem of Neukirch-Uchida for the absolute Galois groups with ramification filtrations of positive characteristic local fields (cf. \cite[Theorem A]{analogue}, \cite[Theorem A]{modified}), and higher local fields with (first) residue characteristics at least $3$ (cf. \cite[Theorems 5, 6]{Ab-fon}).)
This result also shows that, by requiring homomorphisms between the absolute Galois groups of the fields in question to preserve certain structures, we can consider anabelian geometry of fields for which the analogue of the theorem of Neukirch-Uchida fails to hold.
Mochizuki and Hoshi later gave various generalizations of this result (cf. \cite[\S 3]{Topics1}, \cite{note} and \cite[Remark 1.4.1]{intrinsic}).
However, they considered only (fields isomorphic to) $p$-adic local fields.
In the present paper, we consider similar problems for general mixed-characteristic complete discrete valuation fields with perfect residue fields.
(In \cite{GMLF} (resp. \cite{HLF}), the author also considers mono-anabelian reconstruction algorithms of various invariants of these fields (resp. higher local fields) from their absolute Galois groups with ramification filtrations (resp. from their absolute Galois groups).)

One of the main theorems of the present paper, which is a generalization of the results in \cite[\S 3]{Topics1}, is as follows:

\begin{Thm}[{cf. Theorem \ref{main}, Corollary \ref{imp}}]\label{A}
For $i=1,\,2$, let $K_i$ be a mixed-characteristic complete discrete valuation field with perfect residue field, $p_i$ the residue characteristic of $K_i$, $F_i$ the closure of the algebraic closure of $\mathbb{Q}_{p_i}$ in $K_i$, and $G_{K_i}$ (resp. $G_{F_i}$) the absolute Galois group of $K_i$ (resp. $F_i$).
Suppose that $\alpha:G_{K_1}\to G_{K_2}$ is an inertially open (i.e., the image of the inertia subgroup of $G_{K_1}$ is an open subgroup of the inertia subgroup of $G_{K_2}$ (cf. Definition \ref{morphism}, (v))) continuous homomorphism of 01-qLT-type (cf. Definition \ref{morphism}, (i)) (e.g., $\alpha$ is of CHT-type (cf. Definition \ref{morphism}, (iii)), or, more generally, of qLT-type (cf. Definition \ref{morphism}, (i))).

Then there exists a commutative diagram

\[\xymatrix{ G_{K_1} \ar[r]^\alpha \ar@{->>}[d] & G_{K_2}  \ar@{->>}[d] \\ G_{F_1} \ar[r]^{\widetilde{\alpha}} & G_{F_2}}\]
--- where the vertical arrows are the surjections induced by the embeddings $F_i\hookrightarrow K_i$ ($i=1,\,2$).
Moreover, $\widetilde{\alpha}$ is (injective and) geometric (i.e., arises from an embedding $\overline{F_2}\hookrightarrow\overline{F}_1$ that maps to $F_2$ into $F_1$, where $\overline{F_i}$ is an algebraic closure of $F_i$ for $i=1,\,2$ (cf. Definition \ref{morphism}, (iv))).
\end{Thm}

\vskip\baselineskip

Another main theorem of the present paper, which is an analogue of the results in \cite{note}, is as follows:

\begin{Thm}[{cf. Theorem \ref{preserving}}]\label{B}
Let $K_i,\,p_i,\,F_i,\,G_{K_i},\,G_{F_i}$ be as in Theorem \ref{A}, and, for $i=1,\,2$, $F_i^{\mathrm{alg}}$ the algebraic closure of $\mathbb{Q}_{p_i}$ in $K_i$.
Suppose that $\alpha:G_{K_1}\to G_{K_2}$ is an inertially open (cf. Definition \ref{morphism}, (v)) continuous homomorphism of HT-qLT-type (cf. Definition \ref{HT-morphism}, (ii)).
Suppose, moreover, that $F_1^{\mathrm{alg}}$ is Galois over $\mathbb{Q}_p$, where $p:=p_1=p_2$ (cf. the assumption that $\alpha$ is of HT-qLT-type).
Then there exists a(n) (injective and) geometric (cf. Definition \ref{morphism}, (iv)) homomorphism $\alpha':G_{F_1}\to G_{F_2}$.
\end{Thm}

\vskip\baselineskip

Roughly speaking, Theorem \ref{A} treats homomorphisms $\alpha:G_{K_1}\to G_{K_2}$ preserving not only the Hodge-Tate-ness of Galois representations, but also weights appearing in their Hodge-Tate decompositions.
On the other hand, in Theorem \ref{B}, we assume only that $\alpha$ preserves the Hodge-Tate-ness of Galois representations.

\vskip\baselineskip

We shall review the contents of the present paper.
In Section \ref{CHT}, we generalize various notions given in \cite[\S 3]{Topics1} to the general case (i.e., the case of mixed-characteristic complete discrete valuation fields with perfect residue fields), and prove Theorem \ref{A} (cf. Theorem \ref{main}).
Moreover, by applying this result, we prove a certain ``semi-absolute Isom-version'' of the Grothendieck conjecture for hyperbolic curves over mixed-characteristic complete discrete valuation fields with residue fields algebraic over the prime fields (cf. Corollary \ref{semiabs}).
In Section \ref{HT-preserving}, we generalize various notions given in \cite{note} to the general case, and prove Theorem \ref{B} (cf. Theorem \ref{preserving}).
Moreover, we prove a profinite group-theoretic result of the absolute Galois groups of mixed-characteristic complete discrete valuation fields with residue fields algebraic over the prime fields, which is a generalization of \cite[Theorem 2]{JR} (cf. Proposition \ref{hens}).
This result may be of interest independently of the above anabelian results.

\vskip2\baselineskip

\section*{Acknowledgments}

The author would like to thank Professor Akio Tamagawa for numerous helpful advices.
Also, he would like to thank Professors Yuichiro Hoshi, Yuichiro Taguchi and Go Yamashita for giving useful comments, and Professor Atsushi Shiho for giving a question about the author's talk at a conference which suggested him to consider the geometricity of HT-preserving homomorphisms.
The author was supported by the JSPS KAKENHI Grant Number 22J00022.

\vskip2\baselineskip

\section{The geometricity of homomorphisms of CHT, qLT and 01-qLT-type}\label{CHT}

In the present section, we give a generalization of the results in \cite[\S 3]{Topics1}.

\vskip\baselineskip

For a profinite group $G$, we shall say that
\begin{itemize}
\item $G$ is {\it{slim}} if for every open subgroup $H\subset G$, the centralizer $Z_G(H):=\{g\in G\,|\,g\cdot h=h\cdot g,\,\text{for any }h\in H\}$ of $H$ in $G$ is trivial;

\item $G$ is {\it{elastic}} if every topologically finitely generated closed normal subgroup $N\subset H$ of an open subgroup $H\subset G$ of $G$ is either trivial or of finite index in $G$;

\item $G$ is {\it{very elastic}} if $G$ is elastic and not topologically finitely generated.

\end{itemize}

\vskip\baselineskip

Let $K$ be a mixed-characteristic complete discrete valuation field with perfect residue field.
Moreover, we shall write

\begin{itemize}
\item $\overline{K}$ for an algebraic closure of $K$;

\item $\hat{\overline{K}}$ for the completion of $\overline{K}$;

\item $G_K$ for the Galois group $\mathrm{Gal}(\overline{K}/K)$;

\item $I_K\subset G_K$ for the inertia subgroup of $G_K$;

\item $\mathcal{O}_K$ for the ring of integers of $K$;

\item $\mathfrak{M}_K$ for the maximal ideal of $\mathcal{O}_K$;

\item $k=\mathcal{O}_K/\mathfrak{M}_K$ for the residue field of $K$;

\item $p\,(>0)$ for the characteristic of $k$;

\item $\chi_{K}:G_{K}\to\mathbb{Q}_{p}^\times$ for the $p$-adic cyclotomic character.
\end{itemize}

\begin{Def}[{cf. \cite[III, \S A.3\textasciitilde\S A.5]{l-adic}, \cite[Definition 3.1]{Topics1}}]\label{character}

\

\begin{enumerate}
\item[(i)] Let $A$ be an abelian topological group and $\rho,\,\rho':G_K\to A$ continuous homomorphisms.
Then we shall write $\rho\equiv\rho'$ and say that $\rho,\,\rho'$ are {\it{inertially equivalent}} if, for some open subgroup $H\subset I_K$, the restricted homomorphisms $\rho|_H,\,\rho'|_H$ coincide.

\item[(ii)] Let $E$ be a $p$-adic local field all of whose $\mathbb{Q}_p$-conjugates are contained in $K$, $\sigma:E\hookrightarrow K$ an embedding, $\mathcal{O}_E$ the ring of integers of $E$ and $\pi_E$ a uniformizer of $E$.
Define $\chi_{\sigma,\,\pi_E}$ to be the composite homomorphism
\[G_K\to G_E^{\mathrm{ab}}\to \mathcal{O}_E^\times\times\hat{\mathbb{Z}}\twoheadrightarrow \mathcal{O}_E^\times \to \mathcal{O}_E^\times\]
where the first homomorphism is induced by the embedding $\sigma:E\hookrightarrow K$, the second homomorphism is an isomorphism arising from local class field theory (which depends on the choice of $\pi_E$), the third surjection is the projection, and the fourth homomorphism is the inverse automorphism of $\mathcal{O}_E^\times$.
As in \cite[III, \S A.4]{l-adic}, the inertial equivalence class of $\chi_{\sigma,\,\pi_E}$ does not depend on the choice of $\pi_E$.
Thus, we shall often write simply $\chi_{\sigma}$ for $\chi_{\sigma,\,\pi_E}$.

\end{enumerate}

\end{Def}

\vskip\baselineskip

The following proposition is a key ingredient to prove various results in the present paper:

\begin{Prop}[{cf. \cite[III, \S A.5, Corollary]{l-adic}, \cite[Proposition 3.2]{Topics1}}]\label{HT}
Let $E$ be as in Definition \ref{character}, (ii), and $\rho: G_K\to E^\times$ a continuous homomorphism.
Write $V_\rho$ for the $G_K$-module obtained by letting $G_K$ act on $E$ via $\rho$.
Then $\rho$ is Hodge-Tate if and only if
\[\rho\equiv\prod_{\sigma\in\mathrm{Emb}(E,\,K)}\chi_\sigma^{n_\sigma}\]
for some $n_\sigma\in\mathbb{Z}$, where $\mathrm{Emb}(E,\,K)$ is the set of field embeddings $E\hookrightarrow K$ over $\mathbb{Q}_p$.
(For the definition of the Hodge-Tate-ness of Galois representations, see, e.g., \cite[III, \S 1.2, Definition 2]{l-adic}.)
Moreover, in this case, we have an isomorphism of $\hat{\overline{K}}[G_K]$-modules:
\[V_\rho\otimes_{\mathbb{Q}_p}\hat{\overline{K}}\simeq\bigoplus_{\sigma\in\mathrm{Emb}(E,\,K)}\hat{\overline{K}}(n_\sigma),\]
where the ``$(n_\sigma)$'' denotes the $n_\sigma$-th Tate twist.
\end{Prop}

\begin{Rem}
Although \cite[Proposition 3.2]{Topics1} is only stated in the case where $K$ is a $p$-adic local field, a similar result holds in the case where $K$ is a mixed-characteristic complete discrete valuation field with perfect residue field.
\end{Rem}

\vskip\baselineskip

We generalize various notions given in \cite[\S 3]{Topics1} which are defined only for the absolute Galois groups of $p$-adic local fields to those of general mixed-characteristic complete discrete valuation fields (with perfect residue field):

\begin{Def}[{cf. \cite[Definition 3.1]{Topics1}}]
Let $E$ be as in Definition \ref{character}, (ii), and $\rho:G_K\to E^\times$ a continuous homomorphism.

\begin{enumerate}
\item[(i)] We shall say that $\rho$ is {\it{of qLT-type}} (i.e., ``quasi-Lubin-Tate'' type) if there exist an open subgroup $H\subset G_K$, corresponding to a field extension $K_H$ of $K$, and a field embedding $\sigma:E\hookrightarrow K_H$ such that $\rho|_H\equiv\chi_{\sigma}$.

\item[(ii)] We shall say that $\rho$ is {\it{of 01-type}} if it is Hodge-Tate, and, moreover, for any $\sigma\in\mathrm{Emb}(E,\,K)$, $n_\sigma$ defined as in Proposition \ref{HT} belongs to $\{0,\,1\}$.

\item[(iii)] We shall say that $\rho$ is {\it{of ICD-type}} (i.e., ``inertially cyclotomic determinant'' type) if its determinant $\det(\rho):G_K\to\mathbb{Q}_p^\times$ (i.e., the composite of $\rho$ with the norm map $E^\times\to\mathbb{Q}_p^\times$) is inertially equivalent to $\chi_K$.
\end{enumerate}
\end{Def}

\vskip\baselineskip

A similar result to \cite[Proposition 3.3]{Topics1} holds also in the general case:

\begin{Prop}[{cf. \cite[Proposition 3.3]{Topics1}}]\label{qLT-01}
Let $\rho$ and $V_\rho$ be as in Proposition \ref{HT}.
Then the following conditions on $\rho$ are equivalent:

\begin{enumerate}
\item[(i)] $\rho$ is of qLT-type.

\item[(ii)] We have an isomorphism of $\hat{\overline{K}}[G_K]$-modules: $V_\rho\otimes_{\mathbb{Q}_p}\hat{\overline{K}}\simeq\hat{\overline{K}}(1)\oplus\hat{\overline{K}}\oplus\cdots\oplus\hat{\overline{K}}$.

\item[(iii)] $\rho$ is of ICD-type and Hodge-Tate; the resulting $n_\sigma$'s of Proposition \ref{HT} belong to $\{0,\,1\}$.

\item[(iv)] $\rho$ is of ICD-type and of 01-type.
\end{enumerate}
\end{Prop}

\vskip\baselineskip

For $i=1,\,2$, we shall write:

\begin{itemize}
\item $K_i$ for a mixed-characteristic complete discrete valuation field with perfect residue field;

\item $\overline{K_i}$ for an algebraic closure of $K_i$;

\item $\hat{\overline{K_i}}$ for the completion of $\overline{K_i}$;

\item $G_{K_i}$ for the Galois group $\mathrm{Gal}(\overline{K_i}/K_i)$;

\item $I_{K_i}$ for the inertia subgroup of $G_{K_i}$;

\item $k_i$ for the residue field of $K_i$;

\item $p_i\,(>0)$ for the characteristic of $k_i$;

\item $\chi_{K_i}:G_{K_i}\to\mathbb{Q}_{p_i}^\times$ for the $p_i$-adic cyclotomic character;

\item $F_i^{\mathrm{alg}}$ for the algebraic closure of $\mathbb{Q}_{p_i}$ in $K_i$;

\item $F_i$ for the closure of $F_i^{\mathrm{alg}}$ in $K_i$ (note that $F_i$ is a mixed-characteristic complete discrete valuation field with residue field algebraic over $\mathbb{F}_{p_i}$);

\item $\overline{F_i}$ for the algebraic closure of $F_i$ in $\overline{K_i}$;

\item $G_{F_i}$ for the Galois group $\mathrm{Gal}(\overline{F_i}/F_i)$.
\end{itemize}

\begin{Def}[{cf. \cite[Definitions 3.1, (iv) and 3.6, (iii)]{Topics1}}]\label{morphism}
Let $\alpha: G_{K_1}\to G_{K_2}$ be a continuous homomorphism, $L$ the extension of $K_1$ corresponding to $\mathrm{Ker}(\alpha)$ and $\hat{L}$ the completion of $L$.

\begin{enumerate}
\item[(i)] We shall say that $\alpha$ is {\it{of qLT-type}} (resp. {\it{of 01-qLT-type}}) if $p_1=p_2\,(=:p)$, and, moreover, for every pair of open subgroups $H_1\subset G_{K_1}$, $H_2\subset G_{K_2}$ such that $\alpha(H_1)\subset H_2$, and every continuous homomorphism $\rho:H_2\to E^\times$ of qLT-type (where $E$ is a $p$-adic local field all of whose $\mathbb{Q}_{p}$-conjugates are contained in the fields determined by $H_1,\,H_2$), $\rho\circ\alpha|_{H_1}:H_1\to E^\times$ is of qLT-type (resp. of 01-type).

\item[(ii)] We shall say that $\alpha$ is {\it{of HT-type}} (i.e., ``Hodge-Tate'' type) if $p_1=p_2$, and, moreover, the topological $G_{K_1}$-module obtained by composing $\alpha$ with the natural action of $G_{K_2}$ on $\hat{\overline{K_2}}$ is isomorphic (as a topological $G_{K_1}$-module) to $\hat{L}$.

\item[(iii)] We shall say that $\alpha$ is {\it{of CHT-type}} (i.e., ``cyclotomic Hodge-Tate'' type) if $\alpha$ is of HT-type, and,  moreover, the cyclotomic characters of $G_{K_1}$ and $G_{K_2}$ satisfy $\chi_{K_1}=\chi_{K_2}\circ\alpha$.

\item[(iv)] We shall say that $\alpha$ is {\it{geometric}} if it arises from an embedding $\overline{K_2}\hookrightarrow\overline{K_1}$ that maps $K_2$ into $K_1$.

\item[(v)] We shall say that $\alpha$ is {\it{inertially open}} if $\alpha(I_{K_1})$ is an open subgroup of $I_{K_2}$.
\end{enumerate}
\end{Def}

\vskip\baselineskip

The following is the main theorem of the present section:

\begin{Thm}\label{main}
Let $\alpha:G_{K_1}\to G_{K_2}$ be an inertially open continuous homomorphism.
Suppose that $\alpha$ is of qLT-type.
Then there exists a commutative diagram

\[\xymatrix{ G_{K_1} \ar[r]^\alpha \ar@{->>}[d] & G_{K_2}  \ar@{->>}[d] \\ G_{F_1} \ar[r]^{\widetilde{\alpha}} & G_{F_2}}\]
--- where the vertical arrows are the surjections induced by the embeddings $F_i\hookrightarrow K_i$ ($i=1,\,2$).
Moreover, $\widetilde{\alpha}$ is (injective and) geometric.
\end{Thm}

\prf

First, since we have assumed that $\alpha$ is inertially open, the extension $K_2'$ of $K_2$ corresponding to $\alpha(G_{K_1})$ is a henselian discrete valuation field.
So, by replacing $G_{K_2}$ by $\alpha(G_{K_1})$ and $K_2$ by the completion of $K_2'$ if necessary, we may assume that $\alpha$ is surjective.
Since $\alpha$ is of qLT-type, we have $p_1=p_2\,\,(=:p)$.

Let $E$ be a $p$-adic local field, $L_2$ a finite Galois extension of $K_2$ containing all of $\mathbb{Q}_p$-conjugates of $E$, and $L_1$ a finite Galois extension of $K_1$ corresponding to $\alpha^{-1}(G_{L_2})$.
Note that $\alpha$ induces an isomorphism $\overline{\alpha}:\mathrm{Gal}(L_1/K_1)\stackrel{\sim}{\to}\mathrm{Gal}(L_2/K_2)$.
We claim that the following assertion holds:

\begin{quote}
There exists a bijection $\beta_{E,\,L_2}:\mathrm{Emb}(E,\,L_2)\to\mathrm{Emb}(E,\,L_1)$ which is compatible with the actions of $\mathrm{Gal}(L_1/K_1)\simeq\mathrm{Gal}(L_2/K_2)$ (i.e., for any $\gamma\in\mathrm{Gal}(L_1/K_1)=\mathrm{Aut}_{K_1}(L_1)$ and $\sigma_2\in\mathrm{Emb}(E,\,L_2)$, $\beta_{E,\,L_2}(\overline{\alpha}(\gamma)\circ\sigma_2)=\gamma\circ\sigma_1$, where $\sigma_1:=\beta_{E,\,L_2}(\sigma_2)$).
\end{quote}

Indeed, let $L_1'$ be a finite Galois extension of $L_1$ containing all of $\mathbb{Q}_p$-conjugates of $E$, $L_2'$ the finite Galois extension of $L_2$ corresponding to $\alpha(G_{L_1'})$, and $\iota_2:L_2\hookrightarrow L_2'$ the natural inclusion.
Note that there exists a bijection $\mathrm{Emb}(E,\,L_2)\to\mathrm{Emb}(E,\,L_2')$ (given by composing with $\iota_2$).
Fix any $\sigma:E\hookrightarrow L_2$.
Since $\alpha$ is of qLT-type, $\chi_{\iota_2\circ\sigma}\circ\alpha|_{G_{L_1'}}$ is of qLT-type.
So, by enlarging $L_1'$ if necessary, we may assume that there exists a unique embedding $\sigma':E\hookrightarrow L_1'$ such that $\chi_{\iota_2\circ\sigma}\circ\alpha|_{G_{L_1'}}\equiv\chi_{\sigma'}$.
This determines a map $\beta_{E,\,L_2}:\mathrm{Emb}(E,\,L_2)\to\mathrm{Emb}(E,\,L_1')$ (after possibly enlarging $L_1'$ again. Note that $\mathrm{Emb}(E,\,L_2)$ is a finite set).
Since $\mathrm{Gal}(E'/\mathbb{Q}_p)$ acts on $\mathrm{Emb}(E,\,L_1')$ and $\mathrm{Emb}(E,\,L_2)$, and these actions make these sets $\mathrm{Gal}(E'/\mathbb{Q}_p)$-torsors (where $E'$ is the Galois closure of $E/\mathbb{Q}_p$), the map $\beta_{E,\,L_2}$ is bijective (note that, clearly, $\beta_{E,\,L_2}$ is compatible with these actions).
On the other hand, $\alpha$ induces a surjection $\overline{\alpha'}:\mathrm{Gal}(L_1'/L_1)\twoheadrightarrow\mathrm{Gal}(L_2'/L_2)$, and these Galois groups act on $\mathrm{Emb}(E,\,L_1')$ and $\mathrm{Emb}(E,\,L_2)$ (and $\mathrm{Emb}(E,\,L_2')$) respectively.
For any $\gamma'\in\mathrm{Gal}(L_1'/L_1)=\mathrm{Aut}_{L_1}(L_1')$ and any $\sigma\in\mathrm{Emb}(E,\,L_2)$, it is clear that $\chi_{\overline{\alpha'}(\gamma')\circ\iota_2\circ\sigma}\circ\alpha|_{G_{L_1'}}\equiv\chi_{\gamma'\circ\beta_{E,\,L_2}(\sigma)}$.
However, since the image of $\sigma$ is contained in $L_2$, it follows that $\overline{\alpha'}(\gamma')\circ\iota_2\circ\sigma=\iota_2\circ\sigma$ and hence $\gamma\circ\beta_{E,\,L_2}(\sigma)=\beta_{E,\,L_2}(\sigma)$.
Therefore, the image of the embedding $\beta_{E,\,L_2}(\sigma):E\hookrightarrow L_1'$ is contained in $L_1$ for any $\sigma\in\mathrm{Emb}(E,\,L_2)$.
Hence, we obtain a bijection $\beta_{E,\,L_2}:\mathrm{Emb}(E,\,L_2)\to\mathrm{Emb}(E,\,L_1)$.
Moreover, a similar argument shows the compatibility of $\beta_{E,\,L_2}$ with the actions of $\mathrm{Gal}(L_1/K_1)\stackrel{\sim}{\to}\mathrm{Gal}(L_2/K_2)$.
This completes the proof of the above claim.

For $i=1,\,2$, $g\in G_{K_i}$ belongs to $\mathrm{Ker}(G_{K_i}\twoheadrightarrow G_{F_i})$ if and only if the image of $g$ in $\mathrm{Gal}(L_i/K_i)$ acts trivially on $\mathrm{Emb}(E,\,L_i)$ for any possible pair $(E,\,L_2)$.
Thus, $\alpha$ determines an isomorphism $\widetilde{\alpha}:G_{F_1}\stackrel{\sim}{\to}G_{F_2}$.

On the other hand, it is clear that these $\beta_{E,\,L_2}$'s are compatible with each other in the following sense:
\begin{quote}
For any subfield $E'$ of $E$ which is $p$-adic local, and any finite Galois extensions $M_1,\,M_2$ of $K_1,\,K_2$ (corresponding to each other via $\alpha$ and) containing $L_1,\,L_2$, the following diagrams are commutative:
\[\xymatrix{\mathrm{Emb}(E,\,L_2) \ar[r]^{\beta_{E,\,L_2}} \ar[d] & \mathrm{Emb}(E,\,L_1) \ar[d] \\ \mathrm{Emb}(E',\,L_2) \ar[r]^{\beta_{E',\,L_2}} & \mathrm{Emb}(E',\,L_1)}\]
\[\xymatrix{\mathrm{Emb}(E,\,L_2) \ar[r]^{\beta_{E,\,L_2}} \ar[d] & \mathrm{Emb}(E,\,L_1) \ar[d] \\ \mathrm{Emb}(E,\,M_2) \ar[r]^{\beta_{E,\,M_2}} & \mathrm{Emb}(E,\,M_1)}\]
--- where the vertical arrows of the first diagram are determined by the embedding $E'\hookrightarrow E$, and the vertical arrows of the second diagram are determined by the embeddings $L_i\hookrightarrow M_i$ ($i=1,\,2$).
\end{quote}
Therefore, $\alpha$ induces ($\widetilde{\alpha}$ and) a field isomorphism $\alpha_F^{\mathrm{alg}}:F_2^{\mathrm{alg}}\stackrel{\sim}{\to}F_1^{\mathrm{alg}}$ and hence $\alpha_F:F_2\stackrel{\sim}{\to}F_1$.
Moreover, by applying similar arguments to finite extensions of $F_1$ and $F_2$, it follows that the homomorphism between their absolute Galois groups induced by $\alpha_F$ coincides with $\widetilde{\alpha}$.
\qed

\begin{Rem}\label{inertially open}
The assertion obtained by replacing the condition ``$\alpha$ is inertially open'' in the assertion of Theorem \ref{main} by ``$\alpha(G_{K_1})\cap I_{K_2}$ is open in $I_{K_2}$'' also holds.
\end{Rem}

\begin{Rem}
In the situation of Theorem \ref{main}, we have $p_1=p_2\,(=:p)$ (as in the proof of Theorem \ref{main}).
For $i=1,\,2$, the embedding $\mathbb{Q}_p\hookrightarrow K_i$ determines a continuous homomorphism $\phi_i:G_{K_i}\to G_{\mathbb{Q}_p}$ (by considering the algebraic closure of $\mathbb{Q}_p$ in $\overline{K_i}$).
By a similar argument to the proof of Theorem \ref{main}, it follows that $\phi_2\circ\alpha$ coincides with $\phi_1$.
This shows that we can reduce ``absolute situations'' to ``relative situations'' when we consider inertially-open continuous homomorphism of qLT-type between the absolute Galois groups of mixed-characteristic complete discrete valuation fields (with perfect residue fields).
\end{Rem}

\vskip\baselineskip

The following is a generalization of \cite[Theorem 3.5, (i)]{Topics1}:

\begin{Cor}\label{imp}
For $i=1,\,2$, let $K_i$ be a mixed-characteristic complete discrete valuation field with perfect residue field, and $\alpha:G_{K_1}\to G_{K_2}$ an inertially open continuous homomorphism.
Consider the following four conditions:

\begin{enumerate}
\item[(i)] $\alpha$ is geometric.

\item[(ii)] $\alpha$ is of CHT-type.

\item[(iii)] $\alpha$ is of qLT-type.

\item[(iv)] $\alpha$ is of 01-qLT-type.

\end{enumerate}
Then we have the following implications:
\[\text{(i)}\Longrightarrow\text{(ii)}\Longrightarrow\text{(iii)}\Longleftrightarrow\text{(iv)}.\]
If, moreover, for $i=1,\,2$, the residue field of $K_i$ is algebraic over the prime field, the above four conditions are equivalent.
Moreover, in this case, $\alpha$ is injective if one of the above four conditions is satisfied.
\end{Cor}

\prf

The implication (i)$\Longrightarrow$(ii) is clear.
Let us assume that (ii) holds.
Set $p:=p_1=p_2$.
Let $E$ be a $p$-adic local field all of whose $\mathbb{Q}_{p}$-conjugates are contained in $K_2$, and $\rho:G_{K_2}\to E^\times$ a continuous homomorphism of qLT-type.
To prove that (iii) holds, it suffices to prove that $\rho\circ\alpha:G_{K_1}\to E^\times$ is of qLT-type under the condition that all of $\mathbb{Q}_{p}$-conjugates of $E$ are contained in $K_1$.
Since $\alpha$ is of CHT-type, this follows immediately from the equivalence (i)$\Longleftrightarrow$(ii) of Proposition \ref{qLT-01}.
The implication (iii)$\Longrightarrow$(iv) is clear, and the inverse follows from a similar argument to the proof of the implication (b)$\Longrightarrow$(c) of \cite[Theorem 3.5]{Topics1} (note that $\alpha$ is inertially open).
The portion of this corollary concerning the case where the residue field of $K_i$ is algebraic over the prime field for $i=1,\,2$ follows immediately from the above arguments and Theorem \ref{main}.
\qed

\begin{Rem}
Let $\alpha:G_{K_1}\to G_{K_2}$ be an inertially open continuous homomorphism.
Suppose that $\alpha$ is of qLT-type (or, equivalently, of 01-qLT-type (cf. Corollary \ref{imp})).
Since, for $i=1,\,2$, $\chi_{K_i}:G_{K_i}\to\mathbb{Q}_p^\times$ factors through the natural surjection $G_{K_i}\twoheadrightarrow G_{F_i}$ (where $p:=p_1=p_2$), $\alpha$ preserves the $p$-adic cyclotomic characters (i.e., $\chi_{K_1}=\chi_{K_2}\circ\alpha$) by Theorem \ref{main}.
\end{Rem}

\begin{Rem}\label{open}
Suppose that, for $i=1,\,2$, the residue field of $K_i$ is algebraic over the prime field.
Let $\alpha:G_{K_1}\to G_{K_2}$ be a continuous homomorphism such that $\alpha(G_{K_1})\cap I_{K_2}$ is open in $I_{K_2}$ (not necessarily inertially open) (cf. Remark \ref{inertially open}).
By a similar argument to the proof of Theorem \ref{main}, if $\alpha$ satisfies one of (i), (ii) and (iii) in Corollary \ref{imp}, then $\alpha$ is automatically inertially open.
\end{Rem}

\begin{Cor}
For $i=1,\,2$, let $K_i$ be a mixed-characteristic complete discrete valuation field with residue field algebraic over the prime field.
Write $\mathrm{Hom}(K_2,\,K_1)$ for the set of embeddings from $K_2$ to $K_1$, $\mathrm{Hom}^{\mathrm{IO\text{-}CHT}}(G_{K_1},\,G_{K_2})$ (resp. $\mathrm{Hom}^{\mathrm{IO\text{-}qLT}}(G_{K_1},\,G_{K_2})$, resp. $\mathrm{Hom}^{\mathrm{IO\text{-}01\text{-}qLT}}(G_{K_1},\,G_{K_2}))$ for the set of continuous and inertially open homomorphisms of CHT-type (resp. of qLT-type, resp. of 01-qLT-type) from $G_{K_1}$ to $G_{K_2}$, and $\mathrm{Inn}(G_{K_2})$ for the group of inner automorphisms of $G_{K_2}$.
Then the natural maps
\begin{align*}
\mathrm{Hom}(K_2,\,K_1)&\to\mathrm{Hom}^{\mathrm{IO\text{-}CHT}}(G_{K_1},\,G_{K_2})/\mathrm{Inn}(G_{K_2}) \\
&\to\mathrm{Hom}^{\mathrm{IO\text{-}qLT}}(G_{K_1},\,G_{K_2})/\mathrm{Inn}(G_{K_2}) \\
&\to\mathrm{Hom}^{\mathrm{IO\text{-}01\text{-}qLT}}(G_{K_1},\,G_{K_2})/\mathrm{Inn}(G_{K_2})
\end{align*}
are bijective.

\end{Cor}

\prf

The injectivity of the first map follows from \cite[Corollary 4.1]{HoT}.
The surjectivity of the first map and the bijectivity of the second and third maps follow immediately from Corollary \ref{imp}.
\qed

\vskip\baselineskip

In the remainder of this section, for $i=1,\,2$, we shall write

\begin{itemize}
\item $X_i$ for a hyperbolic curve over $K_i$;

\item $\pi_1(X_i)$ for the \'{e}tale fundamental group of $X_i$ (for some choice of basepoint);

\item $\Delta_{X_i}=\pi_1(X_i\times_{\mathrm{Spec}\,K_i}\mathrm{Spec}\,\overline{K_i})$ for the geometric fundamental group of $X_i$ (for some choice of basepoint).
\end{itemize}

\begin{Lem}\label{kernel}
For any isomorphism $\phi:\pi_1(X_1)\to\pi_1(X_2)$ of profinite groups, it holds that $\phi(\Delta_{X_1})=\Delta_{X_2}$.
\end{Lem}

\prf

If $\pi_1(X_1)\,(\simeq\pi_1(X_2))$ is not topologically finitely generated, then $k_1$ and $k_2$ are not finite (cf. \cite[Proposition 1.13]{GMLF}).
(Note that $\Delta_{X_1}$ and $\Delta_{X_2}$ are topologically finitely generated.)
Therefore, both $G_{K_1}$ and $G_{K_2}$ are very elastic (cf. \cite[Proposition 2.1, (iii)]{GMLF}).
So, in this case, the assertion is clear (i.e., $\Delta_{X_i}$ is characterized as the maximal closed normal subgroup of $\pi_1(X_i)$ which is topologically finitely generated for $i=1,\,2$).
On the other hand, the assertion for the case where $\pi_1(X_1)\,(\simeq\pi_1(X_2))$ is topologically finitely generated (hence $k_1$ and $k_2$ are finite) follows from \cite[Lemma 1.3.8]{abs}.
\qed

\begin{Rem}
Let $\phi:\pi_1(X_1)\to\pi_1(X_2)$ be a continuous homomorphism of profinite group.
In the terminology of \cite[Definition 2.4]{Topics1}, similarly to \cite[Proposition 2.5]{Topics1}, $\phi$ is semi-absolute if and only if $\phi$ is pre-semi-absolute (note that $G_{K_2}$ is elastic (cf. \cite[Proposition 2.1, (iii)]{GMLF})).
Moreover, suppose that $k_2$ is not a finite field.
Then $G_{K_2}$ is very elastic (cf. \cite[Propositions 1.13, 2.1, (iii)]{GMLF}).
Therefore, $\phi$ is semi-absolute if and only if $\phi$ is absolute.
\end{Rem}

\vskip\baselineskip

The following is a generalization of \cite[Corollary 3.8]{Topics1} and \cite[Theorem 4.12]{sur}:

\begin{Cor}\label{semiabs}
Suppose that, for $i=1,\,2$, $K_i$ is a mixed-characteristic complete discrete valuation field with residue field algebraic over the prime field.
Write $\mathrm{Isom}(X_1/K_1,\,X_2/K_2)$ for the set of isomorphisms of schemes from $X_1$ to $X_2$ lying over a morphism $\mathrm{Spec}\,K_1\to\mathrm{Spec}\,K_2$, and $\mathrm{Isom}^{\mathrm{IO\text{-}CHT}}(\pi_1(X_1),\,\pi_1(X_2))$ (resp. $\mathrm{Isom}^{\mathrm{IO\text{-}qLT}}(\pi_1(X_1),\,\pi_1(X_2))$, resp. $\mathrm{Isom}^{\mathrm{IO\text{-}01\text{-}qLT}}(\pi_1(X_1),\,\pi_1(X_2))$) for the set of isomorphisms $\phi$ of profinite groups from $\pi_1(X_1)$ to $\pi_1(X_2)$ satisfying the following condition:
\begin{quote}
The isomorphism $\psi:G_{K_1}\to G_{K_2}$ induced by $\phi$ (cf. Lemma \ref{kernel}) is inertially open, and, moreover, of CHT-type (resp. of qLT-type, resp. of 01-qLT-type).
\end{quote}
Moreover, write $\mathrm{Inn}(\pi_1(X_2))$ for the group of inner automorphisms of $\pi_1(X_2)$. Then the natural maps
\begin{align*}
\mathrm{Isom}(X_1/K_1,\,X_2/K_2)&\to\mathrm{Isom}^{\mathrm{IO\text{-}CHT}}(\pi_1(X_1),\,\pi_1(X_2))/\mathrm{Inn}(\pi_1(X_2)) \\
&\to\mathrm{Isom}^{\mathrm{IO\text{-}qLT}}(\pi_1(X_1),\,\pi_1(X_2))/\mathrm{Inn}(\pi_1(X_2)) \\
&\to\mathrm{Isom}^{\mathrm{IO\text{-}01\text{-}qLT}}(\pi_1(X_1),\,\pi_1(X_2))/\mathrm{Inn}(\pi_1(X_2))
\end{align*}
are bijective.
\end{Cor}

\prf

Since $K_i$ is generalized-sub-$p$-adic for $i=1,\,2$ (where $p:=p_1=p_2$), this corollary follows immediately from \cite[Theorem 4.12]{sur} and Corollary \ref{imp}.
\qed

\begin{Rem}
In the situation of Corollary \ref{semiabs}, write $\mathrm{Isom}^{\text{CHT}}(\pi_1(X_1),\,\pi_1(X_2))$ (resp. $\mathrm{Isom}^{\text{qLT}}(\pi_1(X_1),\,\pi_1(X_2))$) for the set of isomorphisms $\phi$ of profinite groups from $\pi_1(X_1)$ to $\pi_1(X_2)$ such that the isomorphism $\psi:G_{K_1}\to G_{K_2}$ is of CHT-type (resp. of qLT-type).
Then, by Remark \ref{open}, the natural inclusions
\begin{align*}
\mathrm{Isom}^{\mathrm{IO\text{-}CHT}}(\pi_1(X_1),\,\pi_1(X_2))&\hookrightarrow\mathrm{Isom}^{\text{CHT}}(\pi_1(X_1),\,\pi_1(X_2)), \\
\mathrm{Isom}^{\mathrm{IO\text{-}qLT}}(\pi_1(X_1),\,\pi_1(X_2))&\hookrightarrow\mathrm{Isom}^{\text{qLT}}(\pi_1(X_1),\,\pi_1(X_2))
\end{align*}
are bijective.
\end{Rem}

\vskip2\baselineskip

\section{A weak anabelian result for HT-preserving homomorphisms}\label{HT-preserving}

In the present section, we consider a similar problem to \cite{note} in a more general situation.
We follow the notations of the previous section.

\begin{Def}[{cf. \cite[Definition 1.3, (i)]{note}}]\label{HT-morphism}
Let $\alpha: G_{K_1}\to G_{K_2}$ be a continuous homomorphism.

\begin{enumerate}
\item[(i)] We shall say that $\alpha$ is {\it{HT-preserving}} (i.e., ``Hodge-Tate-preserving'') if, $p_1=p_2\,(=:p)$, and, moreover, every finite dimensional continuous representation $\phi:G_{K_2}\to\mathrm{GL}_n(\mathbb{Q}_p)$ of $G_{K_2}$ that is Hodge-Tate, the composite $\phi\circ\alpha: G_{K_1}\to\mathrm{GL}_n(\mathbb{Q}_p)$ is Hodge-Tate.

\item[(ii)] We shall say that $\alpha$ is {\it{of HT-qLT-type}} (i.e., ``Hodge-Tate-quasi-Lubin-Tate'' type) (resp. {\it{of weakly HT-qLT-type}} (i.e., ``weakly Hodge-Tate-quasi-Lubin-Tate'' type)) if, $p_1=p_2\,(=:p)$, and, moreover, for every pair of open subgroups $H_1\subset G_{K_1}$, $H_2\subset G_{K_2}$ such that $\alpha(H_1)\subset H_2$, and every $p$-adic local field $E$ all of whose $\mathbb{Q}_{p}$-conjugates are contained in the fields $K_{H_1},\,K_{H_2}$ determined by $H_1,\,H_2$, the composite $\chi_{\sigma_2}\circ\alpha|_{H_1}:H_1\to E^\times$ is Hodge-Tate (resp. is inertially equivalent to a continuous homomorphism $H_1\to E^\times$ that factors through the continuous homomorphism $H_1\to\mathrm{Gal}(\overline{E}/E)$ determined (up to composition with an inner isomorphism) by the embedding $E\stackrel{\sigma_1}{\hookrightarrow}K_{H_1}$), where $\sigma_i:E\hookrightarrow K_{H_i}$ is an embedding for $i=1,\,2$, and $\overline{E}$ is an algebraic closure of $E$.
\end{enumerate}

\end{Def}

\vskip\baselineskip

A similar result to \cite[Lemma 1.4]{note} holds also in the general case:

\begin{Lem}\label{impli}
Let $\alpha:G_{K_1}\to G_{K_2}$ be an inertially open continuous homomorphism.
Consider the following four conditions:

\begin{enumerate}
\item[(i)] $\alpha$ is HT-preserving.

\item[(ii)] For every pair of open subgroups $H_1\subset G_{K_1}$, $H_2\subset G_{K_2}$ such that $\alpha(H_1)\subset H_2$, the restriction $\alpha|_{H_1}:H_1\to H_2$ is HT-preserving.

\item[(iii)] $\alpha$ is of HT-qLT-type.

\item[(iv)] $\alpha$ is of weakly HT-qLT-type.
\end{enumerate}

Then we have the following implications:
\[\text{(i)}\Longleftrightarrow\text{(ii)}\Longrightarrow\text{(iii)}\Longrightarrow\text{(iv)}.\]
\end{Lem}

\begin{Rem}
In the situation of Lemma \ref{impli}, the condition that $\alpha$ is of 01-qLT-type clearly implies that $\alpha$ is HT-preserving (see also Corollary \ref{imp}).

\end{Rem}

\vskip\baselineskip

\begin{Lem}\label{finGal}
Let $K$ be a mixed-characteristic complete discrete valuation field with perfect residue field.
Then $K$ has a finite Galois extension $L$ such that the algebraic closure $F_L^{\mathrm{alg}}$ of $\mathbb{Q}_p$ in $L$ is Galois over $\mathbb{Q}_p$.
\end{Lem}

\prf

Let $\mathbb{F}$ be the algebraic closure of $\mathbb{F}_p$ in $k$ and $F_0$ the quotient field of the Witt ring with coefficients in $\mathbb{F}$.
Then $F^{\mathrm{alg}}$ is a finite extension of $F_0$, where $F^{\mathrm{alg}}$ is the algebraic closure of $\mathbb{Q}_p$ in $K$.
Suppose that $F^{\mathrm{alg}}=F_0(\alpha)$ for $\alpha\in F^\mathrm{alg}$.
Then $\mathbb{Q}_p(\alpha)$ is finite over $\mathbb{Q}_p$ (hence so is the Galois closure $F_0'$ of $\mathbb{Q}_p(\alpha)/\mathbb{Q}_p$).
Then  we can take the finite extension of $K$ determined by the inverse image of the closed normal subgroup $(G_{F^{\mathrm{alg}}}\cap G_{F_0'}=)\,G_{F_0}\cap G_{F_0'}\subset G_{\mathbb{Q}_p}$ by the natural surjection $G_K\twoheadrightarrow G_{F^{\mathrm{alg}}}$ as $L$.
\qed

\vskip\baselineskip

The following is the main theorem of the present section, which is an analogue of \cite[Theorem 3.3]{note} (see also \cite[Remark 1.4.1]{intrinsic}):

\begin{Thm}\label{preserving}
Let $\alpha:G_{K_1}\to G_{K_2}$ be an inertially open continuous homomorphism.
Suppose that $\alpha$ is of HT-qLT-type, and $F_1^{\mathrm{alg}}$ is Galois over $\mathbb{Q}_p$, where $p:=p_1=p_2$.
Then there exists a(n) (injective and) geometric homomorphism $\alpha':G_{F_1}\to G_{F_2}$.
Moreover, for $i=1,\,2$, let $N_i$ be the kernel of the surjection $G_{K_i}\twoheadrightarrow G_{F_i}$ induced by the embedding $F_i\hookrightarrow K_i$.
Then it holds that $\alpha(N_1)\subset N_2$.
\end{Thm}

\prf

We will prove this proposition in a similar way to the proof of \cite[Theorem 3.3]{note}.
First, since we have assumed that $\alpha$ is inertially open, the extension $K_2'$ of $K_2$ corresponding to $\alpha(G_{K_1})$ is a henselian discrete valuation field.
So, by replacing $G_{K_2}$ by $\alpha(G_{K_1})$ and $K_2$ by the completion of $K_2'$ if necessary, we may assume that $\alpha$ is surjective.
For $i=1,\,2$, let us fix an embedding $\sigma_i:\overline{\mathbb{Q}_p}\to \overline{K_i}$ (where $\overline{\mathbb{Q}_p}$ is an algebraic closure of $\mathbb{Q}_p$).
Let $E\subset\overline{\mathbb{Q}_p}$ be a finite Galois extension of $\mathbb{Q}_p$, and $L_1$ and $L_2$ finite Galois extensions of $K_1$ and $K_2$, respectively, such that $\alpha(G_{L_1})=G_{L_2}$, and that $E_i:=\sigma_i(E)\subset L_i$ for $i=1,\,2$.
For $i=1,\,2$, write $\iota_i:E_i\hookrightarrow L_i$ for the natural embedding.
Since $\alpha$ is of HT-qLT-type, the composite
\[\chi_1:G_{L_1}\xrightarrow{\alpha|_{G_{L_1}}}G_{L_2}\xrightarrow{\chi_{\iota_2}}\mathcal{O}_{E_2}^\times\xrightarrow[\sim]{(\sigma_2|_{\mathcal{O}_E^\times})^{-1}}\mathcal{O}_E^\times\xrightarrow[\sim]{\sigma_1|_{\mathcal{O}_E^\times}}\mathcal{O}_{E_1}^\times\]
is Hodge-Tate.
Therefore, by Proposition \ref{HT},
\begin{equation}
\chi_1\equiv\prod_{\tau\in\mathrm{Gal}(E_1/\mathbb{Q}_p)}(\chi_{\iota_1\circ\tau})^{n_\tau}:G_{L_1}\to\mathcal{O}_{E_1}^\times \tag{$\ast$}
\end{equation}
for some $n_\tau\in\mathbb{Z}$ ($\tau\in\mathrm{Gal}(E_1/\mathbb{Q}_p)$).
For $i=1,\,2$, let $J_{K_i}$ (resp. $J_{L_i}$) be the image of $I_{K_i}$ (resp. $I_{L_i}$) in $G_{K_i}^{\mathrm{ab}}$ (resp. $G_{L_i}^{\mathrm{ab}}$), and $\mathrm{Ver}_{L_i/K_i}:G_{K_i}^{\mathrm{ab}}\to G_{L_i}^{\mathrm{ab}}$ the transfer.
Note that $\mathrm{Ver}_{L_i/K_i}$ maps $J_{K_i}$ to $J_{L_i}$.
Therefore, $(\ast)$ shows that there exist open subgroups $S_1\subset J_{K_1}$ and $S_2\subset J_{K_2}$ such that $\alpha^{\mathrm{ab}}(S_1)=S_2$ (where $\alpha^{\mathrm{ab}}:G_{K_1}^{\mathrm{ab}}\to G_{K_2}^{\mathrm{ab}}$ is induced by $\alpha$) and that the following diagram commutes:
\[\xymatrix@C=50pt{S_1 \ar[r] \ar@{->>}[d] & G_{K_1}^{\mathrm{ab}}  \ar[r]^{\mathrm{Ver}_{L_1/K_1}} & G_{L_1}^{\mathrm{ab}} \ar[r]^{\prod_{\tau}(\chi_{\iota_1\circ\tau})^{n_\tau}} & \mathcal{O}_{E_1}^\times \ar[r]^{(\sigma_1|_{\mathcal{O}_E^\times})^{-1}}_\sim & \mathcal{O}_E^\times \ar@{=}[d] \\ S_2 \ar[r] & G_{K_2}^{\mathrm{ab}}  \ar[r]^{\mathrm{Ver}_{L_2/K_2}} & G_{L_2}^{\mathrm{ab}} \ar[r]^{\chi_{\iota_2}} & \mathcal{O}_{E_2}^\times \ar[r]^{(\sigma_2|_{\mathcal{O}_E^\times})^{-1}}_\sim & \mathcal{O}_E^\times.}\]
On the other hand, for $i=1,\,2$, $\mathrm{Gal}(L_i/K_i)$ (resp. $\mathrm{Gal}(E_i/E_i\cap K_i)$) acts on $G_{L_i}^{\mathrm{ab}}$ (by conjugation) (resp. on $\mathcal{O}_{E_i}^\times$).
Note that, for $i=1,\,2$, there exists a natural surjection $\mathrm{Gal}(L_i/K_i)\twoheadrightarrow\mathrm{Gal}(E_iK_i/K_i)=\mathrm{Gal}(E_i/E_i\cap K_i)$, and $\chi_{\iota_2}$ is compatible with these actions.
Moreover, although it is not clear whether $\displaystyle\prod_{\tau}(\chi_{\iota_1\circ\tau})^{n_\tau}$ is compatible with these actions, it is clear that the image of $(G_{L_1}^{\mathrm{ab}})^{\mathrm{Gal}(L_1/K_1)}$ is contained in $(\mathcal{O}_{E_1}^\times)^{\mathrm{Gal}(E_1/E_1\cap K_1)}=\mathcal{O}_{E_1\cap K_1}^\times$.

We claim that, for $i=1,\,2$, the image $T_i$ of $S_i$ in $\mathcal{O}_{E_i}^\times$ is an open subgroup of $\mathcal{O}_{K_i\cap E_i}^\times$.
Indeed, by the definition of the transfer, it is clear that $T_i\subset\mathcal{O}_{K_i\cap E_i}^\times$.
Therefore, to prove the above claim, it suffices to prove that $T_2$ is an open subgroup of $\mathcal{O}_{K_2\cap E_2}^\times$.
On the other hand, $\iota_2:E_2\hookrightarrow L_2$ induces embeddings $\iota_2':E_2\hookrightarrow E_2K_2$ and $\iota_2'':E_2\cap K_2\hookrightarrow K_2$, and these embeddings give open homomorphisms $\chi_{\iota_2'}:G_{E_2K_2}^{\mathrm{ab}}\to\mathcal{O}_{E_2}^\times$ and $\chi_{\iota_2''}:G_{K_2}^{\mathrm{ab}}\to \mathcal{O}_{E_2\cap K_2}^\times$ such that the following diagram commutes:
\[\xymatrix@C=50pt{G_{K_2}^{\mathrm{ab}} \ar[r]^{\mathrm{Ver}_{E_2K_2/K_2}} \ar[d]^{\chi_{\iota_2''}} & G_{E_2K_2}^{\mathrm{ab}} \ar[r]^{\mathrm{Ver}_{L_2/E_2K_2}} \ar[d]^{\chi_{\iota_2'}} & G_{L_2}^{\mathrm{ab}} \ar[d]^{\chi_{\iota_2}} \\ \mathcal{O}_{E_2\cap K_2}^\times \ar[r] & \mathcal{O}_{E_2}^\times \ar[r]^{[L_2:E_2K_2]} & \mathcal{O}_{E_2}^\times}\]
where the first arrow of the second row is the natural inclusion.
Since $\chi_{\iota_2''}(S_2)\subset \mathcal{O}_{E_2\cap K_2}^\times$ is an open subgroup, $T_2$ is also an open subgroup of $\mathcal{O}_{E_2\cap K_2}^\times$, as desired (note that $\mathcal{O}_{E_2\cap K_2}^\times$ is topologically finitely generated).

$\sigma_2\circ\sigma_1^{-1}$ induces an isomorphism $T_1\stackrel{\sim}{\to}T_2$, and this determines a field isomorphism $\beta_{E,\,K_1,\,L_1}(\sigma_1,\,\sigma_2):E_1\cap K_1\stackrel{\sim}{\to} E_2\cap K_2$ (cf. \cite[Lemma 4.1]{version}).
Since we have fixed embeddings $\sigma_1$ and $\sigma_2$, in the remainder of this proof, we shall write simply $\beta_{E,\,K_1,\,L_1}$ for $\beta_{E,\,K_1,\,L_1}(\sigma_1,\,\sigma_2)$.

For various choices of $E$, $K_1$ and $L_1$, it is clear that various $\beta_{E,\,K_1,\,L_1}$'s are compatible in the following sense (note that the restrictions of $\alpha$ to open subgroups of $G_{K_1}$ is also of HT-qLT-type):

\begin{itemize}
\item For any finite Galois extension $K_1'$ of $K_1$ such that the algebraic closure of $\mathbb{Q}_p$ in $K_1'$ is Galois over $\mathbb{Q}_p$, any finite Galois extension $E$ of $\mathbb{Q}_p$, and any finite Galois extensions $L_1$ and $L_2$ of $K_1$ and $K_2$, respectively, such that $\alpha(G_{L_1})=G_{L_2}$, and that $E_i:=\sigma_i(E)\subset L_i$ and $K_i'\subset L_i$ for $i=1,\,2$ (where $K_2'$ is the finite Galois extension of $K_2$ corresponding to $\alpha(G_{K_1'})$), it holds that $\beta_{E,\,K_1',\,L_1}|_{E_1\cap K_1}=\beta_{E,\,K_1,\,L_1}$.

\item For any finite Galois extension $E$ of $\mathbb{Q}_p$, and any finite Galois extensions $L_1,\,L_1'$ and $L_2,\,L_2'$ of $K_1$ and $K_2$, respectively, such that $\alpha(G_{L_1})=G_{L_2}$, $\alpha(G_{L_1'})=G_{L_2'}$, and that $E_i:=\sigma_i(E)\subset L_i$ and $L_i\subset L_i'$ for $i=1,\,2$, it holds that $\beta_{E,\,K_1,\,L_1'}=\beta_{E,\,K_1,\,L_1}$.

\item For any finite Galois extensions $E$ and $E'$ of $\mathbb{Q}_p$ such that $E\subset E'$, and any finite Galois extensions $L_1$ and $L_2$ of $K_1$ and $K_2$, respectively, such that $\alpha(G_{L_1})=G_{L_2}$, and that $E_i':=\sigma_i(E')\subset L_i$ for $i=1,\,2$, it holds that $\beta_{E',\,K_1,\,L_1}|_{E_1\cap K_1}=\beta_{E,\,K_1,\,L_1}$.
\end{itemize}

Therefore, we obtain (a field isomorphism $\beta^{\mathrm{alg}}:F_1^{\mathrm{alg}}\stackrel{\sim}{\to}F_2^{\mathrm{alg}}$ and hence $\beta:F_1\stackrel{\sim}{\to}F_2$, and) a geometric isomorphism $\alpha':G_{F_1}\stackrel{\sim}{\to}G_{F_2}$.
Moreover, the portion of Theorem \ref{preserving} concerning $\alpha(N_1)$ and $N_2$ follows immediately from the above arguments.
\qed

\vskip\baselineskip

By Theorem \ref{preserving}, we obtain a ``weak-Isom'' anabelian result for mixed-characteristic complete discrete valuation fields with residue fields algebraic over the prime fields:

\begin{Cor}
Suppose that, for $i=1,\,2$, $K_i$ is a mixed-characteristic complete discrete valuation field with residue field algebraic over the prime field.
If there exists an inertially open continuous isomorphism $\alpha:G_{K_1}\to G_{K_2}$ of HT-qLT-type, the field $K_1$ is isomorphic to the field $K_2$.
\end{Cor}

\begin{Rem}
Consider the situation of Theorem \ref{preserving} and its proof.
For simplicity, assume that $\alpha$ is surjective.
Then it is not clear whether the following diagram commutes:
\[\xymatrix{ G_{K_1} \ar[r]^\alpha \ar@{->>}[d]^{q_1} & G_{K_2}  \ar@{->>}[d]^{q_2} \\ G_{F_1} \ar[r]^{\alpha'} & G_{F_2}}\]
--- where $q_i$ is the surjection induced by the embedding $F_i\hookrightarrow K_i$ ($i=1,\,2$).
However, it is clear that any open normal subgroup $H\subset G_{K_1}$ such that the algebraic closure of $\mathbb{Q}_p$ in $K_1'$ is Galois over $\mathbb{Q}_p$ (where $K_1'$ is the finite Galois extension of $K_1$ corresponding to $H$) satisfies $\alpha'\circ q_1(H)=q_2\circ\alpha(H)$.
If this equation holds for any open normal subgroup $H\subset G_{K_1}$, then we may prove the following generalization of Theorem \ref{preserving} (without assuming that $F_1^{\mathrm{alg}}$ is Galois over $\mathbb{Q}_p$):

\begin{quote}
Let $\alpha:G_{K_1}\to G_{K_2}$ be an inertially open continuous homomorphism.
Suppose that $\alpha$ is of HT-qLT-type.
Then there exists a commutative diagram
\[\xymatrix{ G_{K_1} \ar[r]^\alpha \ar@{->>}[d]^{q_1} & G_{K_2}  \ar@{->>}[d]^{q_2} \\ G_{F_1} \ar[r]^{\tilde{\alpha}} & G_{F_2}}\]
--- where $q_i$ is the surjection induced by the embedding $F_i\hookrightarrow K_i$ ($i=1,\,2$).
Moreover, $\tilde{\alpha}$ is (injective and) geometric.
\end{quote}
Indeed, let $N_i$ be the kernel of $q_i$ for $i=1,\,2$.
We may assume that $\alpha$ is surjective.
First, consider the case where $F_1^{\mathrm{alg}}$ is Galois over $\mathbb{Q}_p$.
Then we have $\alpha(N_1)=N_2$ by Theorem \ref{preserving}.
Therefore, we obtain an isomorphism $\tilde{\alpha}:G_{F_1}\stackrel{\sim}{\to}G_{F_2}$ such that $\tilde{\alpha}\circ q_1=q_2\circ\alpha$.
Set $\phi:=(\alpha')^{-1}\circ\tilde{\alpha}\in\mathrm{Aut}(G_{F_1})$.
If we have $\alpha'\circ q_1(H)=q_2\circ\alpha(H)$ for any open normal subgroup $H\subset G_{K_1}$, then, for any open normal subgroup $\overline{H}\subset G_{F_1}$, it holds that $\phi(\overline{H})=\overline{H}$.
Therefore, by \cite[Theorem 1]{JR} and Proposition 2.7, $\phi$ is an inner automorphism, and hence geometric.
Since $\alpha'$ is geometric, $\tilde{\alpha}$ is also geometric, as desired.
Next, suppose that $F_1^{\mathrm{alg}}$ is not necessarily Galois over $\mathbb{Q}_p$.
Let $K_1'$ be a finite Galois extension of $K_1$ such that the algebraic closure of $\mathbb{Q}_p$ in $K_1'$ is Galois over $\mathbb{Q}_p$ (cf. Lemma \ref{finGal}), $G_{K_1'}\subset G_{K_1}$ the corresponding open normal subgroup, $K_2'$ the finite Galois extension of $K_2$ corresponding to $\alpha(G_{K_1'})$, and $F_i'$ the closure of the algebraic closure of $\mathbb{Q}_p$ in $K_i'$ for $i=1,\,2$.
We may take $K_1'$ satisfying $\mathrm{Gal}(K_1'/K_1)\simeq\mathrm{Gal}(F_1'/F_1)$ (cf. the proof of Lemma \ref{finGal}).
Note that $\alpha|_{G_{K_1'}}:G_{K_1'}\twoheadrightarrow G_{K_2'}$ is also of HT-qLT-type.
Therefore, we obtain a geometric isomorphism $\tilde{\alpha}':G_{F_1'}\stackrel{\sim}{\to}G_{F_2'}$ which is compatible with $\alpha$ and the natural surjections $q_1':G_{K_1'}\twoheadrightarrow G_{F_1'}$ and $q_2':G_{K_2'}\twoheadrightarrow G_{F_2'}$.
Let $N_i'$ be the kernel of $q_i'$ for $i=1,\,2$.
Since we have assumed that $\mathrm{Gal}(K_1'/K_1)\simeq\mathrm{Gal}(F_1'/F_1)$, $N_1'$ coincides with $N_1$.
We claim that $N_2'$ coincides with $N_2$.
Indeed, clearly we have a surjection $(G_{F_1}\simeq)\,G_{K_2}/N_2'\twoheadrightarrow G_{F_2}$, and this surjection induces a surjection $\mathrm{Gal}(F_1'/F_1)\twoheadrightarrow\mathrm{Gal}(F_2'/F_2)$.
The kernels of these two surjections are isomorphic to each other and hence these kernels are finite groups.
However, since $G_{F_1}$ is slim (cf. \cite[Proposition 2.1, (iii)]{GMLF}), these kernels must be trivial (cf., e.g., \cite[\S 0]{Topics1}).
Therefore, we obtain an isomorphism $\tilde{\alpha}:G_{F_1}\stackrel{\sim}{\to} G_{F_2}$ whose restriction to $G_{F_1'}$ coincides with $\tilde{\alpha}'$, which is geometric.
So, $\tilde{\alpha}$ is also geometric and compatible with $\alpha$, $q_1$ and $q_2$.
\end{Rem}

\begin{Prop}\label{hens}
Let $F$ be a mixed-characteristic henselian discrete valuation field algebraic over $\mathbb{Q}_p$.
Then the absolute Galois group $G_F$ of $F$ is pseudo-$p$-free.
(For the definition of pseudo-$p$-free profinite groups, see \cite[\S 1, Definition]{JR}.)
\end{Prop}

\prf

Let $F_1$ be a finite Galois extension of $F_0:=F$ and $F_2$ a finite Galois extension of $F_1$ which is Galois over $\mathbb{Q}_p$ (note that Lemma \ref{finGal} holds also for mixed-characteristic henselian discrete valuation fields).
Set $F_2=F_0(\alpha)$.
Then $\mathbb{Q}_p(\alpha)$ is finite over $\mathbb{Q}_p$ (hence so is the Galois closure $F'_2$ of $\mathbb{Q}_p(\alpha)/\mathbb{Q}_p$).
Set $F_0':=F_2'\cap F_0$ and denote the intermediate field of $F_2'/F_0'$ corresponding to $\mathrm{Gal}(F_2/F_1)\subset\mathrm{Gal}(F_2/F_0)\simeq\mathrm{Gal}(F_2'/F_0')$ by $F_1'$.
Let $J$ be the image of the wild inertia subgroup $P_{F_1'}\subset G_{F_1'}$ in $G_{F_1'}(p)^{\text{ab-tor}}$ (where, for a profinite group $G$, $G(p)^{\text{ab-tor}}$ denotes the quotient of the maximal abelian pro-$p$ quotient $G(p)^{\mathrm{ab}}$ by the closure of the torsion subgroup).
The inclusion $F_1'\hookrightarrow F_1$ induces a homomorphism $G_{F_1}(p)^{\text{ab-tor}}\to G_{F_1'}(p)^{\text{ab-tor}}$ whose image contains an open subgroup of $J$.
As in the proof of \cite[Theorem 2]{JR}, there exists an element $a\in J$ such that $\{\gamma'(a)\,|\,\gamma'\in\mathrm{Gal}(F_1'/F_0')\}$ (the action of $\mathrm{Gal}(F_1'/F_0')$ on $G_{F_1'}(p)^{\text{ab-tor}}$ is given by conjugation) generates a free $\mathbb{Z}_p$-module of rank $|\mathrm{Gal}(F_1'/F_0')|$.
By multiplying a power of $p$ if necessary, we may assume that $a$ is contained in the image of $G_{F_1}(p)^{\text{ab-tor}}\to G_{F_1'}(p)^{\text{ab-tor}}$, and let $b$ be an element of $G_{F_1}(p)^{\text{ab-tor}}$ which maps to $a$.
Then the closure of the subgroup generated by $\{\gamma(b)\,|\,\gamma\in\mathrm{Gal}(F_1/F_0)\,(\simeq\mathrm{Gal}(F_1'/F_0')\}$ is isomorphic to $\mathbb{Z}_p[\mathrm{Gal}(F_1/F_0)]$.
(Note that the kernel of the natural surjection $G_{F_1}\twoheadrightarrow G_{F_1}(p)^{\text{ab-tor}}$ is a characteristic subgroup of $G_{F_1}$.)
\qed

\vskip2\baselineskip

\end{document}